\documentclass[11pt,a4paper]{article}
\usepackage[nodisplayskipstretch]{setspace}
\usepackage[margin=1in]{geometry}

\usepackage[utf8]{inputenc}
\usepackage{xcolor}
\usepackage{mathtools}
\usepackage{amssymb}
\usepackage{amsthm}
\usepackage{graphicx}
\usepackage{booktabs}
\usepackage[ruled,vlined]{algorithm2e}

\usepackage{mathrsfs}
\usepackage{tikz}
\usepackage{caption}
\usepackage{subcaption}
\usepackage{hyperref}

\newcommand{\B}{\mathbb{B}}
\newcommand{\R}{\mathbb{R}}
\newcommand{\Q}{\mathbb{Q}}

\newcommand{\N}{\mathbb{N}}
\newcommand{\PP}{\mathbb{P}}
\newcommand{\F}{\mathbb{F}}
\newcommand{\E}{\mathbb{E}}
\newcommand{\X}{\mathbb{X}}
\newcommand{\Y}{\mathbb{Y}}
\newcommand{\V}{\mathbb{V}}

\newcommand{\cA}{\mathcal{A}}
\newcommand{\cD}{\mathcal{D}}

\newcommand{\cH}{\mathcal{H}}
\newcommand{\cM}{\mathcal{M}}
\newcommand{\cP}{\mathcal{P}}

\newcommand{\cS}{\mathcal{S}}
\newcommand{\cT}{\mathcal{T}}
\newcommand{\cU}{\mathcal{U}}
\newcommand{\cV}{\mathcal{V}}

\newcommand{\sL}{\mathscr{L}}

\DeclareMathOperator{\dif}{d \!}

\newcommand{\norm}  [1]{\ensuremath{\left  \|       #1  \right \|       }}

\numberwithin{equation}{section}
\theoremstyle{definition}
\newtheorem{defn}{Definition}[section]
\newtheorem{exa}[defn]{Example}
\newtheorem{rem}[defn]{Remark}

\theoremstyle{plain}
\newtheorem{thm}[defn]{Theorem}

\newtheorem{lem}[defn]{Lemma}
\newtheorem{cor}[defn]{Corollary}

\newcommand{\openqed}{\hfill$\square$}

\hypersetup{
    colorlinks,
    linkcolor={red!50!black},
    citecolor={blue!50!black},
    urlcolor={blue!80!black}
}

\title{Multi-objective Markov Decision Process: Set-valued Bellman's principle and fixed point perspective}
\author{Zachary Feinstein\thanks{Stevens Institute of Technology, \href{mailto:zfeinste@stevens.edu}{\texttt{zfeinste@stevens.edu}}}
\and Gabriela Kov{\'{a}}{\v{c}}ov{\'{a}}\thanks{Reykjav{\'{\i}}k University, \href{mailto:gabrielak@ru.is}{\texttt{gabrielak@ru.is}}}
\and Ariel Neufeld\thanks{Nanyang Technological University (NTU), \href{mailto:ariel.neufeld@ntu.edu.sg}{\texttt{ariel.neufeld@ntu.edu.sg}}}
\and Birgit Rudloff\thanks{Vienna University of Economics and Business, \href{mailto:birgit.rudloff@wu.ac.at}{\texttt{birgit.rudloff@wu.ac.at}}}}

\begin{document}

\maketitle
\begin{abstract}
We study infinite-horizon discounted Markov decision processes with multiple cost objectives and finite state and action spaces. We consider the set of value functions generated by randomized Markov policies, where each function describes the expected discounted costs from all initial states under a single policy. We establish a set-valued Bellman equation and characterize this set as the greatest fixed point of the associated Bellman operator. We also show that choosing an attainable cost vector separately for each initial state can produce a function that no single Markov policy realizes. Nevertheless, the attainable costs at individual states satisfy their own Bellman recursion. These results provide a dynamic programming and fixed point framework for multi-objective Markov decision processes.
\end{abstract}

\section{Introduction}
An \emph{infinite-horizon Markov decision problem} (MDP) provides a mathematical framework for sequential decision-making in a stochastic system that operates indefinitely. Let $\mathcal{S}$ denote the state space and $\mathcal{A}$ the action space. At decision time $t=0,1,\ldots$, the system occupies a state $S_t\in\mathcal{S}$, and the decision-maker selects an admissible action $A_t \in \cA(S_t)$. Conditional on $(S_t,A_t)=(s,a)$, the next state $S_{t+1}$ is generated according to the controlled transition probability $P(s'\vert s,a)$, and the transition incurs the immediate cost $c(S_t,A_t,S_{t+1})$. The Markov property requires that, conditional on the current state and action, the distribution of $S_{t+1}$ be independent of the preceding history.

A policy $\pi\in\Pi$ specifies how actions are selected over time, where $\Pi$ denotes the set of admissible policies. Given an initial state $S_0=s$, the expected infinite-horizon discounted cost under policy $\pi$ is
\begin{equation}
	J^{\pi}(s)
	=
	\mathbb{E}^{\pi}\!\left[
	\sum_{t=0}^{\infty}\gamma^{t}
	c(S_t,A_t,S_{t+1}) \,\bigg|\, S_0=s  \right],
\end{equation}
where $\gamma \in (0,1)$ denotes the given discount factor which assigns progressively less weight to costs incurred in the distant future and ensures that the total discounted cost is finite.
The objective is to determine the optimal value $V^*(s):=\inf_{\pi\in\Pi}J^{\pi}(s)$ as well as the corresponding  optimal policy. We refer to \cite{puterman2014markov} and \cite{BR11} for books on (standard, single-objective) Markov decision problems.

Infinite-horizon MDPs arise naturally in continuing operational systems, including inventory and production control, queueing and service operations, maintenance and replacement, transportation, energy management, and healthcare. Brown and Haugh \cite{brown2017information}, for example, study discounted infinite-horizon models for inventory control with autoregressive demand and for dynamic service allocation in a multiclass queue. Financial applications include consumption-investment decisions with a random lifetime, dividend-payment policies for insurance companies, and utility maximization in financial markets; infinite-horizon formulations of these problems are developed by B\"auerle and Rieder \cite[Chapter~9, Sections~9.1--9.3]{BR11}. Credit management provides another example: Liebman \cite{liebman1972markov} formulates the selection of credit-control policies as an infinite-horizon MDP and solves an equivalent linear-programming representation.

The principal solution approach is based on the \emph{dynamic programming principle} (DPP). This allows us to identify the value function $V^*(s)$ as the unique solution (fixed point) of the equation
\begin{equation}
	V^*(s)
	=
	\inf_{a\in\mathcal{A}(s)}
	\sum_{s' \in\mathcal{S}}
	\bigl[c(s,a,s')+\gamma V^*(s')\bigr]
	P(s'\vert s,a),
\end{equation}
known as the \textit{Bellman optimality equation}.
Moreover, the optimal policy can be derived from the optimal action $a^*(s)$ which minimizes the right-hand side of the Bellman equation. This then enables the use of numerical methods such as value iteration or policy iteration to find the optimal value and policy \cite{puterman2014markov}.

In this paper, we study \textit{multi-objective} infinite-horizon Markov decision problems. In contrast to the standard setting discussed above, the decision maker incurs a vector-valued cost $C: \cS \times\cA \times \cS \to \R^n$ in each time period. The goal is to derive the corresponding dynamic programming results, where we pursue an appropriate counterpart of the Bellman's optimality equation for multi-objective MDP and its interpretation in terms of a fixed point.

A number of multi-objective and set-valued dynamic problems have been studied in recent years. These include (finite-horizon) multi-objective optimal control problems, with applications such as mean-risk portfolio selection~\cite{KovacovaRudloff2021}. The validity of dynamic programming is the central question: \cite{KovacovaRudloff2021} establish that a value function of a multi-objective control problem is a \textit{set-valued} object and the counterpart of dynamic programming principle is a \textit{set-valued Bellman's principle}, a recursive relationship for the (set-valued) value function. For continuous-time multi-objective problems, the corresponding counterpart, the set-valued Hamilton-Jacobi-Bellman equations, have been studied~\cite{IseriZhang2023, Visetti2023, HamelVisetti2020}.
 A set-valued Bellman's principle has also been derived in the context of robust multi-objective control~\cite{CialencoKovacova2026},  multivariate and set-valued dynamic risk measures~\cite{FeinsteinRudloff2013, FeinsteinRudloff2014, FeinsteinRudloff2017, AraratFeinstein2021}, Nash equilibria of dynamic games~\cite{FRZ21} and mean-field games~\cite{IseriZhang2024}.

Despite the growing body of work on multi-objective control and dynamic programming, to the best knowledge of the authors, these have not been explored in the context of multi-objective Markov decision processes and infinite horizon. The current work aims to fill this gap by establishing the validity of an appropriate set-valued Bellman's principle and providing an interpretation in terms of a fixed-point framework. By establishing the DPP, this paper provides a solid theoretical framework for future work on reinforcement learning for multi-objective control.

\section{Setting}
\subsection{Preliminaries}
Let $\cT = \{0, 1, \dots \}$ and let $\cS$ be a finite set. The set of probability measures (distributions) over $\cS$ is denoted by $\cM(\cS)$, the power set of $\cS$ is denoted by $\cP(\cS) = 2^\cS$. For a finite space $\cS$ the corresponding $\sigma$-algebra is $\cP(\cS)$, the $\sigma$-algebras of product spaces are obtained as product $\sigma$-algebras. In view of the finiteness assumption, we omit mentions of the $\sigma$-algebras and measurable spaces.

We use the supremum norm $\norm{x}_\infty = \max\limits_{i = 1, \dots, n} |x_i|$ for $x \in \R^n$ and denote by $\B^n_K = \{ x \in \R^n \, \vert \, \norm{x}_\infty \leq K \}$ the closed ball around the origin with radius $K > 0$.
We denote the set of measurable functions $v: \cS \to \R^n$ by $\sL(\cS, \R^n)$, and the notation $\sL(\cS, \B^n_K)$ will be used for the subset of those functions with values in $\B^n_K$. Due to the finiteness of the space $\cS$, the notion of measurability is trivial and, therefore, not addressed.

\subsection{Multi-objective Markov decision process}
Throughout, we consider a state process $(S_t)_{t \in \cT}$ with values in a finite state space $\cS$ and an action (control) process $(A_t)_{t \in \cT}$ with values in a finite action space $\cA$. The transition of the state process fulfills the controlled Markov property with transition probabilities (stochastic kernels) $P: \cS \times \cA \to \cM(\cS)$,
\[
\text{Prob} (S_{t+1} = s' \, \vert \, S_t = s, A_t = a) = P(s' \vert s, a).
\]
Moreover, an initial state distribution $S_0 \sim \nu \in \cM(S)$ is given.

In combination with a scalar-valued cost (or reward) function, we would obtain a standard (single-objective) Markov decision process. In this work, we consider its multi-objective counterpart, a  multi-objective Markov decision process. The cost accrued is multivariate, given by a multi-cost function $C: \cS \times \cA \times \cS \to \R^n$. The aim is to minimize the expected discounted (multi-)cost on an infinite time horizon
\[
\E^\pi \left[ \sum_{t=0}^\infty \gamma^t C(S_t,  A_t,  S_{t+1}) \right],
\]
where $\gamma \in (0,1)$ is a given discount factor. This leads to a multi-objective control problem.

We denote by $\cH_t = \underbrace{(\cS \times \cA) \times \cdots \times (\cS \times \cA)}_{t\text{ times}} \times \cS$ the space of all (state) histories until time $t \in \cT$, and by $\cH = (\cS \times \cA) \times (\cS \times \cA) \times \cdots$ the space of all (state) histories on the infinite time horizon. The (stochastic) process of histories up to time $t$ will be denoted by $H_t = (S_0, A_0, \dots, S_t)$ and its realizations (plays) will be denoted by $h_t = (s_0, a_0, \dots, s_t)$.

\subsection{Multi-objective optimization}
A multi-objective optimization problem is a problem of the form
\begin{align}
\label{MOP}
\tag{MOP}
\text{minimize } F(x) \text{ with respect to } \leq \text{ subject to } x \in \X,
\end{align}
where $\X$ is a non-empty feasible set, $F: \X \to \R^n$ is a vector-valued objective mapping, and $\leq$ is the standard (component-wise) partial order on $\R^n$ induced by the \textit{ordering cone} $\R^n_+$. The \textit{image of the feasible set} of problem~\eqref{MOP} is the set $F[\X] = \{ F(x) \vert x \in \X \}$ and the \textit{upper image} of~\eqref{MOP} is defined as $\cU = \text{cl } (F[\X] + \R^n_+)$.

The prior literature, see e.g.~\cite{KovacovaRudloff2021, CialencoKovacova2026}, has established that a Bellman-type recursion for (finite-horizon) multi-objective control holds for a value function that is a set-valued object corresponding to the image of the feasible set or the upper image. While the upper image is usually preferred due to its set optimization interpretation as an infimum of the multi-objective problem, see~\cite{HHLRS15}, in this work, we opt to utilize the image of the feasible set since boundedness simplifies the fixed point argument. Note that the close relationship between $F[\X]$ and $\cU$ makes it possible to easily extend a recursion for the image of the feasible set into a recursion for the upper image, see e.g.~\cite[Theorem 3.3 or Theorem 4.4]{CialencoKovacova2026}.

A feasible point $\bar x \in \X$ is called a \textit{minimizer} of~\eqref{MOP} if $(F[\bar x] - \R^n_+ \setminus \{0\}) \cap F[\X] = \emptyset$ and it is called a \textit{weak minimizer} if $(F[\bar x] - \text{int } \R^n_+) \cap F[\X] = \emptyset$. A (weak) minimizer is in the literature also referred to as \textit{(weakly) Pareto optimal point} or a \textit{(weakly) efficient point}. The set of images of all (weak) minimizers of \eqref{MOP} is often known as the \textit{(weak) Pareto frontier} or \textit{(weakly) efficient frontier}.

\subsection{Policy classes}

The decision maker controls the MDP through the choice of a policy
$\pi=(\pi_0,\pi_1,\ldots)$. At time $t$, the observed history is
$h_t=(s_0,a_0,s_1,\ldots,a_{t-1},s_t)$.
A general history-dependent randomized decision rule is a stochastic kernel
$\pi_t(\cdot\mid h_t)$ on the admissible action set $\mathcal A(s_t)$. We distinguish
between \emph{history-dependent} policies, for which the decision rule $\pi_t$ may
depend on the full history $h_t$, and \emph{Markov} policies, for which $\pi_t$ depends
only on the current state $s_t$. We also distinguish between \emph{deterministic}
policies, where the decision rule selects a single admissible action, and
\emph{randomized} policies, where the decision rule specifies a probability distribution
over admissible actions. Thus, a randomized Markov decision rule is of the form
$\pi_t(\cdot\mid s_t)$, whereas a deterministic Markov decision rule is induced by a
measurable selector $A_t$ satisfying $A_t(s_t)\in \cA(s_t)$.

For single-objective discounted MDPs, restricting attention to \textit{Markov} policies is
without loss of optimality under the standard assumptions of infinite-horizon MDP
theory. In particular, the optimal value obtained over history-dependent policies
coincides with the optimal value obtained over Markov policies
\cite[Chapter~7]{BR11}. Moreover, if for every state the infimum in the corresponding Bellman
equation is attained by an admissible action and a measurable minimizing selector can be
chosen, then an optimal deterministic Markov policy exists. In finite-state and
finite-action discounted MDPs, this attainment condition is automatic; in this setting,
there exists an optimal deterministic stationary policy \cite[Theorem~6.2.7]{puterman2014markov}.
For the multi-objective setting considered in this paper, \cite[Theorem~3]{CMH}
provides the corresponding reduction at the level of Pareto images: the image generated
by Pareto-optimal history-dependent policies coincides with the image generated by
Pareto-optimal Markov policies. Hence, restricting attention to Markov policies does not
remove any Pareto-optimal objective vector.

The role of \textit{randomization} is different. In a single-objective  MDP with an
expected discounted-cost criterion, randomization does not improve the optimal value.
More precisely,
 optimizing the corresponding Bellman operator over randomized decision
rules yields the same value as optimizing over deterministic decision rules
\cite[Proposition~6.2.1]{puterman2014markov}. Intuitively, for a fixed state, the expected
one-step cost plus continuation value under a randomized action is a convex combination
of the corresponding quantities under deterministic actions; hence no randomized action
can improve upon the best deterministic action. Consequently, if the infimum in the Bellman equation is
attained by an admissible action and a measurable minimizer can be selected, an optimal
randomized Markov policy can be replaced by a deterministic Markov policy with the same
scalar value. If no minimizing action exists, this is a nonattainment issue for the
scalar optimization problem rather than a strict improvement obtained through
randomization.

In multi-objective MDPs, however, the object of interest is not a single scalar value
but the attainable set of a vector-valued criteria. Since the underlying order is only
partial, randomization can generate objective vectors that are not attainable by
deterministic policies, for instance through convex combinations of deterministic
objective vectors. Such randomized objective vectors may be Pareto optimal and therefore
relevant for the efficient frontier, as illustrated in Example~\ref{ex:1}. This is a
genuine distinction from the scalar case: the issue is not merely whether an infimum is
attained, but whether all Pareto-relevant objective vectors are attainable within the
chosen policy class. Since deterministic policies may fail to generate the full
Pareto-relevant image, we work throughout this paper with randomized Markov policies.

\begin{exa}{\cite[Example~1]{CMH}} \label{ex:1}
Consider the multi-objective MDP given in Figure~\ref{fig:1} with initial state $S_0 = s_0 = B$ and multi-cost function $C(s, a, s') = C(s)$. In this case, the expected discounted multi-cost is fully determined by the action $A_0$. Only two expected discounted costs can be achieved across all deterministic policies, $\left( \frac{\gamma}{1 - \gamma}, 0 \right)$ for $\pi_0(B) = l$ and $\left( 0, \frac{\gamma}{1 - \gamma} \right)$ for $\pi_0(B) = r$. However, consider a randomized policy where each action $l$ and $r$ is chosen with equal probability in stage $t=0$. This is a Pareto optimal policy and its expected cost $\left( \frac{\gamma}{2(1 - \gamma)}, \frac{\gamma}{2(1 - \gamma)} \right)$ cannot be achieved by any deterministic policy. \openqed
\end{exa}
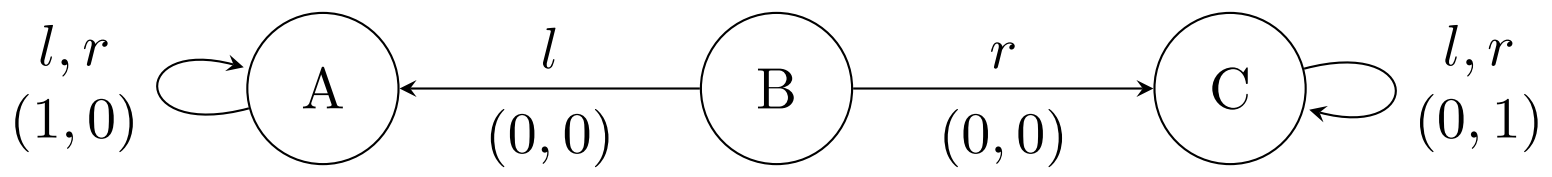
\begin{figure}  \begin{center}
\caption{\label{fig:1} Multi-objective MDP of Example~\ref{ex:1}}
\begin{tikzpicture}[->, >=stealth, node distance=3cm]

\node[circle, draw, minimum size=1cm] (A) {A};
\node[circle, draw, minimum size=1cm, right of=A] (B) {B};
\node[circle, draw, minimum size=1cm, right of=B] (C) {C};

\draw (B) --
    node[midway, above] {$l$}
    node[midway, below] {$(0,0)$}
(A);

\draw (B) --
    node[midway, above] {$r$}
    node[midway, below] {$(0,0)$}
(C);

\draw (A) edge[loop left] node[left, align=center] {$l, r$\\ $(1,0)$} (A);
\draw (C) edge[loop right] node[right, align=center] {$l, r$\\ $(0,1)$} (C);

\end{tikzpicture}
\end{center} \end{figure}

In the remainder of this subsection, we formalize admissible policies.
Let $\cA(\cdot): \cS \rightrightarrows \cA$ be a set-valued mapping with non-empty values, where $\cA(s)$ represents the admissible actions for state $s \in \cS$. We denote by
\[
\cD = \left\lbrace \varpi: \cS \to \cM(\cA) \quad \vert \quad \forall s \in \cS: \; \varpi(\cA(s) \vert s)  = 1 \right\rbrace
\]
the set of admissible randomized (Markov) decision rules and by
\[
\Pi = \left\lbrace \pi = (\pi_0, \pi_1, \dots ) \quad \vert \quad \forall t \in \cT: \;  \pi_t \in \cD \right\rbrace
\]
the set of all admissible randomized Markov policies. For a policy $\pi \in \Pi$ and a decision rule $\varpi \in \cD$ we define a forward shift operator $\theta(\pi) = (\pi_1, \pi_2, \dots)$ and a concatenation $(\varpi, \pi) = (\varpi, \pi_0, \pi_1, \dots)$.

Observe that these operators preserve admissibility. Namely, for any policy $\pi \in \Pi$ and decision rule $\varpi \in \cD$, it holds
\begin{align}\label{eq:policies}
    \left\lbrace \theta(\pi) \; \vert \; \pi \in \Pi \right\rbrace = \Pi = \left\lbrace (\varpi, \pi) \; \vert \; \varpi \in \cD, \; \pi \in \Pi \right\rbrace.
\end{align}

\subsection{Probability measure}
We follow \cite[Section 2.2]{HL12} to establish the probability measure used by the decision maker who implements policy $\pi$. Since the results are standard, we provide them without proofs.

By the Ionescu-Tulcea theorem \cite[Proposition C.10]{HL12}, for each policy $\pi \in \Pi$  there exists a unique probability measure $\PP^\pi_\nu$ on $\cH = (\cS \times \cA) \times (\cS \times \cA) \times \cdots$ such that
\begin{align*}
&\PP^\pi_\nu (S_0 = s) = \nu(s), \\
&\PP^\pi_\nu (A_t = a \vert H_t = h_t) = \pi_t(a \vert s_t), \\
&\PP^\pi_\nu (S_{t+1} = s \vert H_t = h_t, A_t = a_t) = P(s \vert s_t, a_t),
\end{align*}
where $h_t = (s_0, a_0, \dots, s_t)$. The expectation under $\PP^\pi_\nu$ will be denoted $\E^\pi_\nu$. In case of a Dirac measure $\nu = \delta_s$, we will use the notation $\E^\pi [ \,\cdot\, \vert S_0 = s]$.

With a slight abuse of notation, we use $\theta: \cH \to \cH$ to denote a forward shift operator given by $\theta ((s_0, a_0, s_1, a_1, \dots)) = (s_1, a_1, s_2, a_2, \dots)$ as well as a forward shift operator for the policies. In the following lemma we summarize the decomposition and Markov type properties on the level of probability measures.

\begin{lem} \label{lemma:prob}
Let $\PP^\pi_\nu$ be  the probability measure on $\cH$ induced by policy $\pi \in \Pi$, prior $\nu \in \cM(\cS)$ and transition kernel $P$. The following holds for the distribution of $H= (S_0, A_0, S_1, \dots)$:
\vspace{-\topsep}
\begin{enumerate}
\item[(i)] The marginal law of $(S_0, A_0, S_1)$ is given by
\[
\PP^\pi_\nu (S_0 = s_0, A_0 = a_0, S_1 = s_1) = \nu(s_0) \pi_0(a_0 \vert s_0) P(s_1 \vert s_0, a_0)
\]
for all $(s_0, a_0, s_1) \in \cH_1$.
\item[(ii)] The conditional tail distribution of $\theta(H) \vert S_1 = s$ is given by
\[
\PP^\pi_\nu\big( \theta(H) \in \cdot \,\big|\, S_1 = s \big) = \PP^{\theta(\pi)}_{\delta_s}(\cdot).
\]
\item[(iii)] For any $h_1 = (s_0, a_0, s_1) \in \cH_1$ and any $B \subseteq \cH$, a decomposition is given by
\[
\PP^\pi_\nu  \big( H_1 = h_1, \theta(H) \in B  \big) = \PP^\pi_\nu  \big( H_1 = h_1 \big) \PP^{\theta(\pi)}_{\delta_{s_1}}(B).
\]

\item[(iv)] Given a function $f: \cS \times (\cA \times \cS)^k \to \R^n$ it holds
\begin{align*}
&\E^\pi \left[ f(S_0, A_0, \dots, S_k) \, \vert \, S_0 = s_0 \right]  \\
&= \int_{(\cA \times \cS)^k} f(s_0, a_0, \dots, s_k) P(\dif s_k \vert s_{k-1}, a_{k-1}) \pi_{k-1} (\dif a_{k-1} \vert s_{k-1}) \dots P(\dif s_1 \vert s_0, a_0) \pi_0 (\dif a_0 \vert s_0).
\end{align*}
\end{enumerate}
\end{lem}

\section{From policies to values}
Throughout this section, let the admissible policy $\pi \in \Pi$ be fixed. We introduce the expected discounted multi-cost associated with policy $\pi$ and derive its recursive form, which will be used to derive a Bellman-type recursion for the multi-objective problem.

For a policy $\pi$ we define a policy (state) value function $V^\pi: \cS \to \R^n$ as
\[
V^\pi (s) :=  \E^\pi \left[ \sum_{t=0}^\infty \gamma^t C(S_t,  A_t,  S_{t+1}) \; \Big| \; S_0 = s\right] \quad \quad \text{ for all } s \in \cS.
\]
The policy state value function $V^\pi$ satisfies a recursive relationship known as the \textit{Bellman consistency equation}. Observe that this recursive form for multi-objective MDP is a direct counterpart to the single-objective case.
\begin{thm} \label{thm:piBellman}
Given an admissible policy $\pi \in \Pi$, for all states $s \in \cS$ it holds
\begin{align*}
V^\pi (s) &= \E^\pi \left[ C(S_0, A_0, S_1) + \gamma V^{\theta(\pi)} (S_1) \; \mid \; S_0 = s \right] \\
&= \sum_{a \in \cA} \pi_0(a \vert s) \sum_{s' \in \cS} P(s' \vert s, a) \left( C(s, a, s') + \gamma V^{\theta(\pi)} (s') \right).
\end{align*}
\end{thm}
Before proving the theorem, we first need the following lemma.
\begin{lem} \label{lemma:VQ}
Define a function $Q^\pi: \cS \times \cA \to \R^n$ for all $(s,a) \in \cS \times \cA$ with $\pi_0(a|s) > 0$ as
\[
Q^\pi (s, a) :=  \E^\pi \left[ \sum_{t=0}^\infty \gamma^t C(S_t,  A_t,  S_{t+1}) \; \Big| \; S_0 = s, A_0 = a \right],
\]
and set $Q^\pi (s, a) := 0$ otherwise.
For the state value function $V^\pi$ and the state-action value function $Q^\pi$ it holds
\begin{align} \label{eq:1}
V^\pi (s) = \sum_{a \in \cA} \pi_0(a \vert s) Q^\pi (s, a)
\end{align}
and for $(s,a)\in \cS \times \cA$ satisfying $\pi_0(a|s) > 0$ it holds
\begin{align} \label{eq:2}
Q^\pi (s, a) = \sum_{s' \in \cS} P(s' \vert s, a) \left( C(s, a, s') + \gamma V^{\theta(\pi)} (s') \right).
\end{align}
\end{lem}
\begin{proof}
The tower property of the expectation is used throughout this proof. The absolute integrability of the series follows from the finiteness of the spaces $\cS, \cA$.
We use Lemma~\ref{lemma:prob}. For \eqref{eq:1}, let $s \in \cS$ and observe that
\begin{align*}
V^\pi (s) &=  \E^\pi \left[ \sum_{t=0}^\infty \gamma^t C(S_t,  A_t,  S_{t+1}) \; \Big| \; S_0 = s\right] =  \E^\pi \left[ \E^\pi \left[ \sum_{t=0}^\infty \gamma^t C(S_t,  A_t,  S_{t+1}) \; \Big| \; S_0, A_0\right] \Big| \; S_0 = s\right] \\
&= \E^\pi \left[ Q^\pi (S_0,A_0) \; \vert \; S_0 = s\right] = \sum_{a \in \cA} \pi_0(a \vert s) Q^\pi (s, a).
\end{align*}
To prove \eqref{eq:2}, observe that, in view of Lemma~\ref{lemma:prob}(ii), it holds
\[
\E^{\pi} \left[ \sum_{t=1}^\infty \gamma^{t-1} C(S_t,  A_t,  S_{t+1}) \; \vert \; S_1 = s\right] = \E^{\theta(\pi)} \left[ \sum_{t=0}^\infty \gamma^t C(S_t,  A_t,  S_{t+1}) \; \vert \; S_0 = s\right] = V^{\theta(\pi)} (s).
\]
Then by the Markov property, for $(s,a) \in \cS \times \cA$ with $\pi_0(a|s) > 0$ we obtain
\begin{align*}
Q^\pi (s, a) &=  \E^\pi \left[ \sum_{t=0}^\infty \gamma^t C(S_t,  A_t,  S_{t+1}) \; \Big| \; S_0 = s, A_0 = a \right] \\
&=  \E^\pi \left[ C(S_0, A_0, S_1) + \gamma \sum_{t=1}^\infty \gamma^{t-1} C(S_t,  A_t,  S_{t+1}) \; \Big| \; S_0 = s, A_0 = a \right] \\
&=  \E^\pi \left[ C(S_0, A_0, S_1) + \gamma \E^\pi \left[ \sum_{t=1}^\infty \gamma^{t-1} C(S_t,  A_t,  S_{t+1}) \; \Big| \; S_1\right] \; \Big| \; S_0 = s, A_0 = a \right] \\
&=  \E^\pi \left[ C(S_0, A_0, S_1) + \gamma V^{\theta(\pi)} (S_1) \; \Big| \; S_0 = s, A_0 = a \right] \\
&= \sum_{s' \in \cS} P(s' \vert s, a) \left( C(s, a, s') + \gamma V^{\theta(\pi)} (s') \right).
\end{align*}
\end{proof}

\begin{proof}[Proof of Theorem~\ref{thm:piBellman}]
In view of \eqref{eq:1} and \eqref{eq:2} from Lemma~\ref{lemma:VQ}, we obtain
\begin{align*}
V^\pi (s) &= \sum_{a \in \cA} \pi_0(a \vert s) Q^\pi (s, a) = \sum_{a \in \cA} \pi_0(a \vert s) \sum_{s' \in \cS} P(s' \vert s, a) \left( C(s, a, s') + \gamma V^{\theta(\pi)} (s') \right),
\end{align*}
respectively
\begin{align*}
V^\pi (s) &= \E^\pi \left[ Q^\pi (S_0,A_0) \; \vert \; S_0 = s\right] = \E^\pi \left[ \E^\pi \left[ C(S_0, A_0, S_1) + \gamma V^{\theta(\pi)} (S_1) \; \Big| \; S_0, A_0 \right] \; \vert \; S_0 = s\right] \\
&= \E^\pi \left[  C(S_0, A_0, S_1) + \gamma V^{\theta(\pi)} (S_1) \; \Big| \; S_0 = s\right].
\end{align*}
\end{proof}

Before moving to the optimal value function, we provide bounds on the value associated with an arbitrary policy.
Given the finite state and action spaces, there exists a finite bound
$$ \bar{C} := \max\limits_{s, s' \in \cS, a \in \cA}  \norm{C(s, a, s')}_\infty < \infty $$
for the supremum norm (as well as each component) of the multi-cost function. We show that the value associated with any policy is bounded (in the supremum norm) by $M := \frac{\bar{C}}{1-\gamma}$.

\begin{lem}
\label{lemma:bound}
For an admissible policy $\pi \in \Pi$, it holds $V^\pi \in \sL(\cS, \B^n_M)$, i.e.,
\begin{align*}
\norm{V^\pi (s)}_\infty &\leq \frac{\bar{C}}{1-\gamma} \quad \text{ for all states } s \in \cS.
\end{align*}
\end{lem}
\begin{proof}
For an arbitrary history $(S_0, A_0, S_1, A_1, \dots)$, it holds
\begin{align*}
\norm{\sum_{t=0}^\infty \gamma^t C(S_t,  A_t,  S_{t+1})}_\infty &\leq \max_{i = 1, \dots, n} \sum_{t=0}^\infty \gamma^t \left\vert C_i(S_t,  A_t,  S_{t+1}) \right\vert \leq  \frac{\bar{C}}{1-\gamma}.
\end{align*}
The result follows.
\end{proof}

\begin{rem}
\begin{itemize}
\item Bounds on the multi-cost would be similarly obtained if the finiteness of the state and action spaces was replaced by a compactness assumption alongside mild continuity of the cost function.

\item For the state-action value function $Q^\pi$ we obtain an identical bound $\norm{Q^\pi (s, a)}_\infty \leq \frac{\bar{C}}{1-\gamma}$ across all $(s,a) \in \cS \times \cA$, that is, $Q^\pi \in \sL(\cS \times \cA, \B^n_M)$.

\end{itemize}
\end{rem}

\section{Optimal value function}
We now return to the multi-objective problem faced by the decision maker, the problem of minimizing the expected discounted multi-cost across all admissible Markov policies,
\begin{align} \label{problem}
\min_{\pi \in \Pi}  \quad \; V^\pi (S_0) =  \E^\pi \left[ \sum_{t=0}^\infty \gamma^t C(S_t,  A_t,  S_{t+1}) \; \vert \; S_0 \right] \quad \text{ w.r.t. } \leq_{\R^n_+}.
\end{align}
The aim of this work is to establish a Bellman-type recursive relationship for this multi-objective control problem. It has been established in the literature in the context of finite-horizon problems that the counterpart of the dynamic programming principle for \textit{multi-objective} controls is a \textit{set-valued Bellman's principle}, see e.g.~\cite{KovacovaRudloff2021}.

We define the (optimal) value function of problem~\eqref{problem} as the set (collection) of functions providing the policy values
\[
\V^* = \left\lbrace s \mapsto V^\pi(s) \; \mid \; \pi \in \Pi \right\rbrace \subseteq \sL(\cS, \R^n)
\]
across all admissible policies. This section delivers the main result in two steps. First, we prove the validity of the \textit{set-valued Bellman optimality equation} for the value function of the infinite horizon problem. Second, we interpret this optimality equation in terms of a \textit{fixed point}. These results show that the structure of dynamic programming results for infinite horizon problems is preserved also in the context of multiple objectives, when the value function is interpreted as a set of functions.

\begin{rem}
In the context of multi-objective control, the optimal value function is often treated as a \textit{set-valued function} attaining the upper image of the multi-objective problem, see~\cite{KovacovaRudloff2021, CialencoKovacova2026}. In this paper, we instead work with \textit{a set of functions} -- specifically, the collection of policy-valued functions across all admissible policies -- which enables an easy interpretation of the set-valued Bellman's principle in terms of a fixed point.  See Appendix~\ref{appendix} for a further discussion.
\end{rem}

\subsection{Set-valued Bellman optimality equation}
For single-objective control on infinite horizon, the dynamic programming principle has the form of a \textit{Bellman optimality equation}, which shows that the optimal value function can be obtained through one-step recursions. We show that this remains true also in the multi-objective case, when the optimal value function is treated as a set of functions.

\begin{thm} \label{thm:Bellman}
It holds
\begin{align} \label{eq:Bellman}
\V^* &= \left\lbrace \; s \longmapsto \sum_{a \in \cA} \varpi(a \vert s) \sum_{s' \in \cS} P(s' \vert s, a) \left( C(s, a, s') + \gamma V (s') \right) \; \Big| \; \varpi \in \cD, \, V \in \V^* \; \right\rbrace.
\end{align}
\end{thm}

The proof of Theorem~\ref{thm:Bellman} will be given after Lemma~\ref{lemma:VQopt}.
\begin{lem} \label{lemma:VQopt}
Let us define the set $\Q^*$ of state-action values as
\[
\Q^* := \left\lbrace (s, a) \mapsto Q^\pi(s, a) \; \mid \; \pi \in \Pi \right\rbrace.
\]
Then, it holds
\begin{align} \label{eq:3}
\Q^* = \left\lbrace  \; (s, a) \mapsto  \sum_{s' \in \cS} P(s' \vert s, a) \left( C(s, a, s') + \gamma V (s') \right) \; \Big| \;  V \in \V^* \; \right\rbrace
\end{align}
and
\begin{align} \label{eq:4}
\V^* = \left\lbrace \; s \mapsto \sum_{a \in \cA} \varpi(a \vert s) Q(s,a) \; \Big| \; \varpi \in \cD, \, Q \in \Q^* \; \right\rbrace.
\end{align}
\end{lem}
\begin{proof}
In view of~\eqref{eq:policies}, it holds $\V^* = \{ V^\pi \; \vert \; \pi \in \Pi\} = \{ V^{\theta(\pi)} \; \vert \; \pi \in \Pi\}$. Then, Lemma~\ref{lemma:VQ} yields
\begin{align*}
\Q^* &= \left\lbrace (s, a) \mapsto  Q^\pi(s, a) \; \mid \; \pi \in \Pi \right\rbrace \overset{\eqref{eq:2}}{=} \left\lbrace  (s, a) \mapsto  \sum_{s' \in \cS} P(s' \vert s, a) \left( C(s, a, s') + \gamma V^{\theta(\pi)} (s') \right) \; \mid \; \pi \in \Pi \right\rbrace \\
&= \left\lbrace (s, a) \mapsto  \sum_{s' \in \cS} P(s' \vert s, a) \left( C(s, a, s') + \gamma V (s') \right) \; \mid \; V \in \V^* \right\rbrace.
\end{align*}
This proves \eqref{eq:3}. For \eqref{eq:4}, we use Lemma~\ref{lemma:VQ} to observe
\begin{align*}
\V^* &= \left\lbrace s \mapsto V^\pi(s) \; \mid \; \pi \in \Pi \right\rbrace \overset{\eqref{eq:1}}{=}  \left\lbrace s \mapsto  \sum_{a \in \cA} \pi_0(a \vert s) Q^\pi (s, a) \; \mid \; \pi \in \Pi \right\rbrace \\ &\subseteq \left\lbrace \; s \mapsto  \sum_{a \in \cA} \varpi(a \vert s) Q(s,a) \; \Big| \; \varpi \in \cD, \, Q \in \Q^* \; \right\rbrace.
\end{align*}
For the reverse inclusion, let $\varpi \in \cD$ and $Q \in \Q^*$. By definition of the collection $\Q^*$, there exists a policy $\pi \in \Pi$ such that $Q = Q^\pi$. In view of~\eqref{eq:policies}, also policy $(\varpi, \theta(\pi)) \in \Pi$ is admissible and since $Q^\pi$ does not depend on the decision rule $\pi_0$, it follows that $Q^\pi = Q^{(\varpi, \theta(\pi))}$. Finally, the above together with Lemma~\ref{lemma:VQ} implies
\begin{align*}
&\left\lbrace \; s \mapsto \sum_{a \in \cA} \varpi(a \vert s) Q(s,a)  \;  \Big| \; \varpi \in \cD, \, Q \in \Q^* \; \right\rbrace = \left\lbrace \; s \mapsto  \sum_{a \in \cA} \varpi(a \vert s) Q^\pi(s,a) \; \Big| \; \varpi \in \cD, \, \pi \in \Pi \; \right\rbrace \\
&\quad = \left\lbrace \; s \mapsto \sum_{a \in \cA} \varpi(a \vert s) Q^{(\varpi, \theta(\pi))}(s,a) \; \Big| \; \varpi \in \cD, \, \pi \in \Pi \; \right\rbrace \\
&\quad \overset{\eqref{eq:1}}{=}  \left\lbrace \; s \mapsto   V^{(\varpi, \theta(\pi))}(s) \; \Big| \; \varpi \in \cD, \, \pi \in \Pi \; \right\rbrace \subseteq \V^*.
\end{align*}
\end{proof}

\begin{proof}[Proof of Theorem~\ref{thm:Bellman}]
The set-valued Bellman optimality equation follows from Lemma~\ref{lemma:VQopt},
\begin{align*}
\V^* &\overset{\eqref{eq:4}}{=} \left\lbrace \; s \mapsto  \sum_{a \in \cA} \varpi(a \vert s) Q(s,a) \; \Big| \; \varpi \in \cD, \, Q \in \Q^* \; \right\rbrace \\ &\overset{\eqref{eq:3}}{=} \left\lbrace \; s \mapsto  \sum_{a \in \cA} \varpi(a \vert s) \sum_{s' \in \cS} P(s' \vert s, a) \left( C(s, a, s') + \gamma V (s') \right) \; \Big| \; \varpi \in \cD, \, V \in \V^* \; \right\rbrace.
\end{align*}
\end{proof}

In the next subsection, we will interpret the set-valued Bellman optimality equation in terms of a \textit{fixed point}. For that purpose, we define a mapping $f: \cD \times \sL(\cS, \R^n) \to \sL(\cS, \R^n)$ as
\begin{align} \label{eq:mapf}
f(\varpi, v) (s) = \int_{\cA} \int_{\cS}  \left( C(s, a, s') + \gamma v (s') \right) P(\dif s' \vert s, a) \varpi(\dif a \vert s).
\end{align}
We note that the recursion \eqref{eq:Bellman} can be expressed in terms of the mapping $f$ as
\[
\V^* = \left\lbrace \; f(\varpi, v) \; \vert \; \varpi \in \cD, \, v \in \V^* \; \right\rbrace.
\]
In Lemma~\ref{lemma:bound} we have shown that a policy value $V^\pi$ corresponding to an admissible policy $\pi$ is bounded, i.e., $V^\pi \in \sL(\cS, \B^n_M)$. Now we show that the mapping $f$ preserves this bound.

\begin{lem}\label{lemma:boundF}
If $\varpi \in \cD$ and $v \in \sL(\cS, \B^n_M)$, then $f(\varpi, v) \in \sL(\cS, \B^n_M)$.
\end{lem}
\begin{proof}
For $v \in \sL(\cS, \B^n_M)$ it holds $|v_i(s')| \leq M$ for all $s' \in \cS$ and all $i = 1, \dots, n$. Then,
\begin{align*}
\norm{f(\varpi, v) (s)}_\infty = &\norm{\sum_{a \in \cA} \varpi(a \vert s) \sum_{s' \in \cS} P(s' \vert s, a) \left( C(s, a, s') + \gamma v (s') \right)}_\infty \\
&\leq \; \max_{i=1, \dots, n} \sum_{a \in \cA} \varpi(a \vert s) \sum_{s' \in \cS} P(s' \vert s, a) \cdot \big| C_i(s, a, s') + \gamma v_i (s') \big| \\
&\leq  \sum_{a \in \cA} \varpi(a \vert s) \sum_{s' \in \cS} P(s' \vert s, a) \cdot (\bar{C} + \gamma M) = M
\end{align*}
\end{proof}

\subsection{Fixed point}
\label{scFP}
For single-objective infinite horizon control problems, the Bellman optimality equation allows an interpretation as a \textit{fixed point}. We show that the same is true also for the set-valued counterpart in case of a multi-objective infinite horizon control problem.

In light of the set-valued Bellman optimality equation \eqref{eq:Bellman} and Lemma~\ref{lemma:boundF}, we define a mapping $\F: \cP(\sL(\cS, \B^n_M)) \to \cP(\sL(\cS, \B^n_M))$ as
\begin{align*}
\F(\V) :=& \left\lbrace \; s \mapsto \sum_{a \in \cA} \varpi(a \vert s) \sum_{s' \in \cS} P(s' \vert s, a) \left( C(s, a, s') + \gamma V (s') \right) \; \Big| \; \varpi \in \cD, \, V \in \V \; \right\rbrace \\
=& \left\lbrace \; f(\varpi, v) \; \Big| \; \varpi \in \cD, \, v \in \V \; \right\rbrace.
\end{align*}
Observe that $\F$ maps a set of functions onto a set of functions and Lemma~\ref{lemma:boundF} guarantees that if $\V \subseteq \sL(\cS, \B^n_M)$, then $\F(\V) \subseteq \sL(\cS, \B^n_M)$.

The set-valued Bellman optimality equation \eqref{eq:Bellman} stated in terms of the mapping $\F$ corresponds to $\V^* = \F(\V^*)$. The \textit{Tarski fixed point framework} allows us to interpret the value function $\V^*$ as a fixed point of the mapping $\F$. In contrast to the scalar setting, herein we do not obtain uniqueness of the fixed point. However, we prove below that the optimal value function $\V^*$ corresponds to the \textit{largest} fixed point of the mapping $\F$.

\begin{cor}\label{tarski}
Let $g: \X \times \Y \rightarrow \X$ and define $G:\cP(\X) \rightarrow \cP(\X)$ by $G(X) := \{g(x,y) \, \vert \, x \in X, y \in Y_0\}$ for any $X \in \cP(\X)$ and some fixed $Y_0 \subseteq \Y$.
The set of fixed points of $G$ (i.e., $\{X^* \subseteq \X \; | \; X^* = G(X^*)\}$) forms a complete lattice under the subset order relation $\subseteq$.
In particular, there exists a largest fixed point $\bar X = G(\bar X)$ such that $X^* \subseteq \bar X$ for any $X^* = G(X^*)$.
\end{cor}
\begin{proof}
Note that $(\cP(\X),\subseteq)$ forms a complete lattice where the infimum is defined by the intersection and the supremum by the union (as the largest subset of a collection of sets is the intersection and the smallest superset is the union).
Additionally, trivially, $G(X) \subseteq G(X')$ for any $X \subseteq X' \in \cP(\X)$.
The result trivially follows by the Tarski fixed point theorem \cite[Theorem 1.11]{AB06}.
\end{proof}
\begin{rem}
The smallest fixed points for $G$ as constructed in Corollary~\ref{tarski} is, trivially, the empty set $\emptyset$.
\end{rem}

In order to prove that the optimal value function $\V^*$ is the largest fixed point of the mapping $\F$, we rewrite the mapping as $\F(\V) = \{f(\varpi, v) \, \vert \, \varpi \in \cD, v \in \V\}$.

\begin{thm}
\label{largestFix}
Let $\V \subseteq \sL(\cS, \B^n_M)$ be a fixed point of $\F$, i.e., it holds $\V=\F(\V)$.
Then $\V\subseteq \V^*$. Hence, $\V^*$ is the largest fixed point of $\F$.
\end{thm}
\begin{proof}
Consider an arbitrary $v \in \V$, where $\V$ is a fixed point of $\F$. Then, there exists $\varpi_0 \in \cD$ and $v^{(1)} \in \V$ such that $v = f(\varpi_0, v^{(1)})$. By the same argument, $v^{(1)} = f(\varpi_1, v^{(2)})$ for some $\varpi_1 \in \cD$ and $v^{(2)} \in \V$. We can repeat this argument $k$ times for any $k \in \N$ and obtain
\begin{align*}
v(s_0) &= f(\varpi_0, f(\varpi_1, \dots f(\varpi_{k-1}, v^{(k)}))) (s_0) \\
&= \int_{\cA \times \cS} \dots \int_{\cA \times \cS} \left( \sum_{t = 0}^{k-1} \gamma^t C(s_t, a_t, s_{t+1}) + \gamma^{k} v^{(k)} (s_k) \right) P(\dif s_k \vert s_{k-1}, a_{k-1}) \varpi_{k-1}(\dif a_{k-1} \vert s_{k-1}) \\
& \hspace{7cm} \dots P(\dif s_1 \vert s_0, a_0) \varpi_0(\dif a_0 \vert s_0)
\end{align*}
for all $s_0 \in \cS$. In view of Lemma~\ref{lemma:prob}(iv), it holds
\begin{align*}
v(s_0) &= \E^\pi \left[ \sum_{t = 0}^{k-1} \gamma^t C(S_t, A_t, S_{t+1}) + \gamma^{k} v^{(k)} (S_k) \; \Big| \; S_0 = s_0 \right]
\end{align*}
for a policy $\pi \in \Pi$ satisfying $(\pi_0, \dots, \pi_{k-1}) = (\varpi_0, \dots, \varpi_{k-1})$. In particular, we obtain that for $v \in \V$, there exists a policy $\pi \in \Pi$, such that for all $s_0 \in \cS$, it holds
\begin{align} \label{k-iter}
v(s_0) &= \lim\limits_{k \to \infty} \E^\pi \left[ \sum_{t = 0}^{k-1} \gamma^t C(S_t, A_t, S_{t+1}) + \gamma^{k} v^{(k)} (S_k) \; \Big| \; S_0 = s_0 \right].
\end{align}

Now, let $|\mathcal{S}| \in \N$ denote the cardinality of the finite state space $\mathcal{S}$ and let $\mathcal{S}=\{s_1,\dots, s_{|\mathcal{S}|}\}$ denote its elements. Then, for each $k \in \N$, we can define an $n|\mathcal{S}|$-dimensional vector $\mathfrak{v}^k:=(v^{(k)}(s_1),\dots, v^{(k)}(s_{|\mathcal{S}|}))\in \R^{n|\mathcal{S}|}$.
 Note that the sequence of vectors $(\mathfrak{v}^k)_{k \in \N}$ is uniformly bounded since each element of the sequence $ (v^{(k)})_{k \in \N}\subseteq \V$ is uniformly bounded by $M$. Therefore, by Heine-Borel, there exists a subsequence $(\mathfrak{v}^{k_j})_{j \in \N}$, and an element $ \mathfrak{v}_\infty \in \R^{n|\mathcal{S}|}$ such that $\lim_j  \mathfrak{v}^{k_j} = \mathfrak{v}_\infty$. In particular, this implies for every $j\in \N$
 and $s\in \mathcal{S}$
that
\begin{equation*}
	\norm{v^{(k_j)}(s)}_\infty\leq \norm{\mathfrak{v}^{k_j}}_\infty\leq \sup_\ell \norm{\mathfrak{v}^{k_\ell}}_\infty<\infty.
\end{equation*}
This, together with dominated convergence and the assumption that $0<\gamma <1$, ensures that
\[
\lim_{j \to \infty} \E^\pi \left[ \gamma^{k_j} v^{(k_j)} (S_{k_j}) \; \vert \; S_0 = s_0 \right] = 0.
\]
Therefore, by \eqref{k-iter} and using that the multi-cost function $C$ is uniformly bounded by $\bar{C}<\infty$, we obtain by dominated convergence that
\begin{align*}
v(s_0) &= \lim_{j\rightarrow\infty} \E^\pi \left[ \sum_{t = 0}^{k_j-1} \gamma^t C(S_t, A_t, S_{t+1}) + \gamma^{k_j} v^{(k_j)} (S_{k_j}) \; \vert \; S_0 = s_0 \right] \\
&= \E^\pi \left[ \sum_{t = 0}^{\infty} \gamma^t C(S_t, A_t, S_{t+1}) \; \vert \; S_0 = s_0 \right] = V^\pi (s_0).
\end{align*}
Hence, $v = V^\pi \in \V^*$ and $\V\subseteq \V^*$. Since $\V$ was an arbitrary fixed point of $\F$, it follows that $\V^*$ is the largest set being a fixed point of $\F$ and contains all other fixed points of $\F$, which means by the subset relation that $\V^*$ is the largest fixed point of $\F$.
\end{proof}

\appendix
\section{Optimal value function: Set of functions versus set-valued function}
\label{appendix}
We define the optimal value function as the \textit{set of all policy value functions}
\[
\V^* = \left\lbrace s \mapsto V^\pi(s) \; \mid \; \pi \in \Pi \right\rbrace \subseteq \sL(\cS, \R^n).
\]
In contrast, other works (e.g.~\cite{KovacovaRudloff2021, CialencoKovacova2026}) define the optimal value function of a multi-objective control problem as a \textit{set-valued function}. In the context of the present problem, such (alternative) optimal value function would be $\cV^*: \cS \rightrightarrows \R^n$ given by
\[
\cV^* (s) = \left\lbrace V^\pi(s) \; \mid \; \pi \in \Pi \right\rbrace = \left\lbrace v(s) \; \mid \; v \in \V^* \right\rbrace.
\]
Clearly, all elements of $\V^*$ are \textit{selectors} of the set-valued function $\cV^*$, i.e.\ $\V^* \subseteq \text{Sel } \cV^*$, where $\text{Sel } \cV^*$ denotes the set of all selectors of $\cV^*$. However, Example~\ref{ex2} illustrates that not all selectors of $\cV^*$ belong to $\V^*$, i.e.\ not all selectors of $\cV^*$ are attainable by a single admissible policy across all states. Thus, $\V^* \subsetneq \text{Sel } \cV^*$ can hold. Following, e.g., \cite[Theorem 2.1.6]{molchanov2005theory}, this can only happen if $\V^*$ is not a \textit{decomposable} set, i.e., there exists a pasting $s \mapsto V(s) := V^{\pi}(s)\mathbf{1}_{\{s \in D\}} + V^{\bar\pi}(s)\mathbf{1}_{\{s \in D^c\}}$ for some $\pi,\bar\pi \in \Pi$ and $D \subseteq \cS$ such that $V \not\in \V^*$. We leave the study of when we recover decomposability of $\V^*$ for future work. Here $\mathbf{1}_{\{ \cdot\}}$ denotes the indicator function.

\begin{figure}[h]
\begin{center}
\begin{tikzpicture}[->, >=stealth, node distance=3cm]

\node[circle, draw, minimum size=1cm] (A) {A};
\node[circle, draw, minimum size=1cm, right of=A] (B) {B};

\draw (A) edge[bend left=25] node[midway, above] {0.5} (B);
\draw (B) edge[bend left=25] node[midway, below] {0.5} (A);

\draw (A) edge[loop left] node[left] {0.5} (A);
\draw (B) edge[loop right] node[right] {0.5} (B);

\end{tikzpicture}
\end{center}
\caption{ \label{fig:2} Transition probabilities of MDP in Example~\ref{ex2}.}
\end{figure}
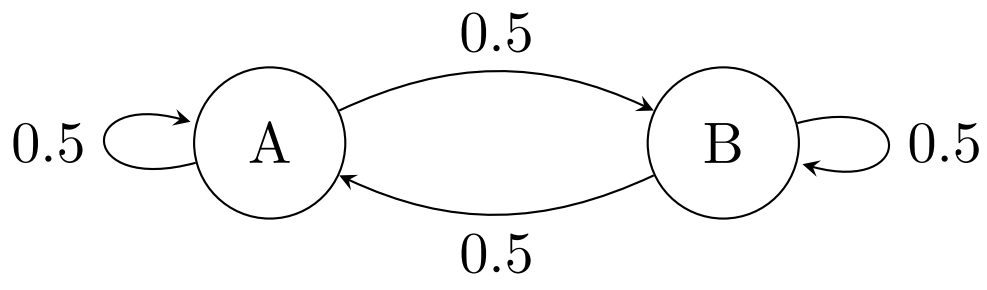
\begin{exa} \label{ex2}
 Consider a simple MDP with two states $\cS = \{A, B\}$ and two actions $\cA = \{l, r\}$. The transition probabilities are independent of the chosen action and are depicted in Figure~\ref{fig:2}, the cost function depends only on the action and is given by
$C(s,a,s') = C(a) =
\bigl(\begin{smallmatrix} 1 \\ 0 \end{smallmatrix}\bigr) \mathbf{1}_{\{a = l\}}
+ \bigl(\begin{smallmatrix} 0 \\ 1 \end{smallmatrix}\bigr) \mathbf{1}_{\{a = r\}}.$

The policy $\pi^l$ of always choosing action $l$ yields the value
$V^{\pi^l} (A) = V^{\pi^l} (B)
= \sum_{t=0}^\infty \gamma^t
\bigl(\begin{smallmatrix} 1 \\ 0 \end{smallmatrix}\bigr)
= \bigl(\begin{smallmatrix} \frac{1}{1-\gamma} \\ 0 \end{smallmatrix}\bigr).$
 Similarly, the policy $\pi^r$ of always choosing action $r$ yields
$V^{\pi^r} (A) = V^{\pi^r} (B) = \bigl(\begin{smallmatrix} 0 \\ \frac{1}{1-\gamma}  \end{smallmatrix}\bigr)$.
 Then, a function $v: \cS \to \R^2$ given by
\[
v(A) = \begin{pmatrix}
\frac{1}{1-\gamma} \\ 0
\end{pmatrix}, \quad \quad v(B) = \begin{pmatrix}
0 \\\frac{1}{1-\gamma}
\end{pmatrix}
\]
is a \textit{selector} of the set-valued function $\cV^*$, i.e., $v(s) \in \cV^*(s)$ for all $s \in \cS$. But $v$ is not produced by any admissible policy, i.e., $v \not\in \V^*$.
\end{exa}

The main body of this paper derived the appropriate dynamic programming principle for the value function $\V^*$ in Theorem~\ref{thm:Bellman}. A natural follow up question is whether also $\cV^*$ admits a recursive form.  As a consequence  of the Bellman optimality equation \eqref{eq:Bellman} for $\V^*$, for all states $s \in \cS$, it holds
\begin{align} \label{eq:recursion1}
\cV^* (s) = \left\lbrace \; \sum_{a \in \cA} \varpi(a \vert s) \sum_{s' \in \cS} P(s' \vert s, a) \left( C(s, a, s') + \gamma V (s') \right) \; \Big| \; \varpi \in \cD, \, V \in \V^* \; \right\rbrace.
\end{align}
Note, however, that \eqref{eq:recursion1} is not a full recursion for the set-valued function $\cV^*$ since it involves only a subset $\V^* \subsetneq \text{Sel } \cV^*$ of the selectors of $ \cV^*$.
However, a true recursion, and thus a set-valued Bellman's principle for the set-valued function $\cV^*$
can be obtained.
\begin{thm}\label{thm:recApp}
For all states $s \in \cS$, it holds
  \begin{align}\label{eq:recursion2}
\cV^* (s) = \left\lbrace \; \sum_{a \in \cA} \varpi(a \vert s) \sum_{s' \in \cS} P(s' \vert s, a) \left( C(s, a, s') + \gamma V (s') \right) \; \Big| \; \varpi \in \cD, \, V \in \text{Sel } \cV^* \; \right\rbrace.
\end{align}
\end{thm}
We provide a \textit{sketch of the proof}. Consider a pair $\varpi \in \cD, \, V \in \text{Sel } \cV^*$ generating an element of the right-hand side of equation~\eqref{eq:recursion2}. Given a selector, for every state $s_1 \in \cS$ there exists a (Markov randomized) policy $\pi^{(s_1)} \in \Pi$ such that $V(s_1) = V^{\pi^{(s_1)}} (s_1)$. Consequently, the considered right-hand side element of equation~\eqref{eq:recursion2} is generated as an expectation under a \textit{history-dependent} randomized policy $\bar \pi$ given by $\bar \pi_0(\cdot \vert s_0) = \varpi (\cdot \vert s_0)$ and $\bar \pi_t (\cdot \vert s_1, s_t) = \pi^{(s_1)}_{t-1} (\cdot \vert s_t)$ for $t > 0$.

To show that the right-hand side element belongs to the set $\cV^* (s)$, it needs to be shown that an objective value arising from a history-dependent randomized policy can also be obtained with a Markov randomized policy. \cite[Theorem 3]{CMH} provides this result in a setting with multi-cost (more precisely, reward) functions $\tilde C: \cS \times \cA \to \R^n$. Defining $\tilde{C} (s, a) = \int_{\cS} C(s,a, s') P(\dif s' \vert s,a)$ for all pairs $(s,a)$ reduces our problem to the setting of~\cite{CMH}.  We omit the remaining details for the sake of brevity. \qed

Note that it is quite common in the multi-objective case to obtain a set-valued Bellman's principle for several objects related to the multi-objective problem. For instance, in \cite{KovacovaRudloff2021} is was shown that for the mean-risk problem the \textit{upper images} recurs backwards in time, but also the \textit{image sets} as well as the \textit{weakly efficient frontiers}. The focus in \cite{KovacovaRudloff2021} was then on the upper image as it is a convex set, and thus the most convenient to work with and to approximate and the best choice for a numerical implementation of the recursion.

The above Bellman's principle \eqref{eq:recursion2} for the set-valued function $\cV^*$ corresponds to the recursion of the image sets. However, as the focus in the present paper is on deriving for the first time a Bellman's principle for a multi-objective infinite-horizon Markov decision problem and provide in total analogy to the scalar case an interpretation of it as a fixed point, the set-valued Bellman optimality equation \eqref{eq:Bellman} for the set $\V^*$ of all policy value functions is the right choice, as it allows for a direct application of Tarski fixed point theorem in Section~\ref{scFP}. However, for future work it might be the case that other formulations of the Bellman's principle for a multi-objective infinite-horizon Markov decision problem like the one in terms of the set-valued function $\cV^*$ derived in Theorem~\ref{thm:recApp} might be relevant as well.

\bibliographystyle{alpha}
\bibliography{biblio}

\end{document}